\documentclass[12pt,oneside,reqno]{amsart}

\usepackage{amsmath}
\usepackage{amsthm}

\usepackage{aliascnt}
\usepackage{amssymb}
\usepackage{latexsym}
\usepackage{amsfonts,setspace}
\usepackage{fullpage}
\usepackage{blkarray}
\usepackage{indentfirst}
\usepackage{algorithmic}
\usepackage{algorithm}
\usepackage{color}
\usepackage{changepage}
\usepackage{paralist}
\usepackage{bbold}
\usepackage{mathtools}
\usepackage{mathabx,epsfig}
\usepackage{tikz}
\usepackage{bm}
\usepackage{verbatim}
\let\oldtextit\textit 
\renewcommand\emph[1]{\oldtextit{\color{RoyalBlue}#1}}
\definecolor{RoyalBlue}{cmyk}{1, 0.50, 0, 0}
\usepackage{todonotes}

\usepackage{makecell}
\usepackage{pict2e}
\usepackage{hyperref}
\usepackage{cleveref}
\crefname{prop}{proposition}{propositions}
    \crefname{lem}{lemma}{Lemmas}

    \theoremstyle{definition}
    \newtheorem{defi}{Definition}
    
\newtheorem{thm}{Theorem}    
    \newaliascnt{lem}{thm}
    \newtheorem{lem}[lem]{Lemma}
    \aliascntresetthe{lem}
    
    \newaliascnt{rmk}{thm}
    \newtheorem{rmk}[rmk]{Remark}
    \aliascntresetthe{rmk}
    
    \newaliascnt{prop}{thm}
    \newtheorem{prop}[prop]{Proposition}
    \aliascntresetthe{prop}
    
    \newaliascnt{conj}{thm}
    
    \aliascntresetthe{conj}
    \newaliascnt{example}{thm}
    \newtheorem{example}[example]{Example}
    \aliascntresetthe{example}
    
    \newaliascnt{cor}{thm}
    \newtheorem{cor}[cor]{Corollary}
    \aliascntresetthe{cor}

\newcommand{\RR}{\mathbb{R}}

\newcommand{\ZZ}{\mathbb{Z}}

\newcommand{\mainfilecheck}[1]{0}

\newcommand{\spanz}{\text{span}}

\normalbaselines

\title{Uniquely balanced polyhedral complexes are\\ highly connected}
\DeclareMathOperator{\relint}{relint}
\DeclareMathOperator{\rec}{rec}

\newcommand{\Tp}{\mathcal T(p)}

\author{Kisun Lee}
\address{Department of Mathematical Sciences, Korea Advanced Institute of Science and Technology}
\email{kisunl@kaist.ac.kr}
\urladdr{https://klee669.github.io}

\date{}

\begin{document}

\begin{abstract}
We establish higher connectivity of polyhedral complexes admitting a unique positive balanced top cycle. Our main result shows that, after quotienting by common lineality, every skeleton satisfies the corresponding higher connectivity property. The proof bounds the number of connected components surviving facet deletions by the dimension of the top-cycle space using homological and linear-algebraic arguments. As applications, we prove higher connectivity for all skeleta of arbitrary pointed polyhedral structures on tropical linear spaces of valuated matroids and Bergman fans including nonrealizable cases. These results answer questions raised by Maclagan and Yu.
\end{abstract}
\maketitle

\section{Introduction}

The tropicalization of an irreducible algebraic variety is connected through codimension one \cite{cartwright2013connectivity}: any two maximal cells can be joined by a sequence of maximal cells meeting along codimension-one faces. This connectivity underlies algorithms for computing tropical varieties \cite{bogart2007computing}. \cite{maclagan2022higher} strengthens this property by proving higher connectivity: for a $d$-dimensional irreducible variety over a field of characteristic zero, its tropicalization is $(d-\ell)$-connected through codimension one where $\ell$ is the dimension of the lineality space of the chosen polyhedral complex (see \cite{gandini2021toric} for positive characteristic). This means that connectivity persists after deleting any $d-\ell-1$ closed maximal cells. Its significance lies in the use of algebraic geometry methods, notably a tropical Bertini theorem, to establish strong combinatorial consequences from algebraic irreducibility.

\cite[Section 5]{maclagan2022higher} asks whether tropical irreducibility suffices for higher connectivity. We approach this question through the \emph{space of balanced top-dimensional cycles}. For a pure $d$-dimensional polyhedral complex $\Delta$, let $\ker\partial_d^\Delta$ denote its top-cycle space, defined using the cellular boundary map with coefficients in the top exterior powers of the direction spaces of its facets. We say that $\Delta$ is \emph{uniquely balanced} if it admits a positive top cycle and $\dim\ker\partial_d^\Delta=1$, as its positive top cycle is unique up to scaling. Our first main result shows that this condition forces higher connectivity, without any rationality assumption.

\begin{thm}\label{thm:general}
Let $\Delta$ be a finite pure $d$-dimensional real polyhedral complex admitting a positive top cycle, and suppose $\dim\ker\partial_d^\Delta=1$. Then, the lineality spaces of all cells of $\Delta$ coincide, defining a subspace $U$ such that $\Delta/U$ is pointed. For every $1\leq k\leq d-\dim U$, its $k$-skeleton is $k$-connected through codimension one. Equivalently, the $(\dim U+k)$-skeleton of $\Delta$ is $k$-connected through codimension one.
\end{thm}

In the usual rational setting of tropical geometry, the uniquely balanced condition is equivalent to tropical irreducibility. Consequently, our theorem establishes higher connectivity for tropically irreducible balanced complexes. The distinction from algebraic irreducibility is essential because an algebraically irreducible variety may not have a tropically irreducible tropicalization. For example, the irreducible curve
\[X=V((x+1)(y+1)+t)\subset (\mathbb C((t))-\{0\})^2\]
has tropicalization $\mathrm{Trop}(X)=\{u=0\}\cup\{v=0\}$ where $(u,v)$ are the tropical coordinates. This tropicalization has two independent balanced top cycles. Conversely, the Bergman fan of the V\'amos matroid is uniquely balanced but cannot arise as the tropicalization of an algebraic variety. It illustrates that our work applies beyond the scope of the algebraic irreducibility result (see \Cref{ex:vamos}).

Whereas \cite{maclagan2022higher} uses a tropical Bertini theorem, we prove by filling boundaries on deleted facets and recovering the balancing condition through linear algebra. The key ingredient is \Cref{prop:components} which bounds the number of connected components surviving the deletion of at most $d-1$ facets by the dimension of the top-cycle space. The idea is to construct an independent top cycle from each surviving component by filling the boundary of the restricted positive cycle on the deleted facets. We apply \Cref{thm:general} to tropical linear spaces of valuated matroids and Bergman fans. In particular, we establish higher connectivity for all skeleta of arbitrary pointed polyhedral structures on these spaces, including nonrealizable cases (see \Cref{thm:main,cor:bergmanConnectivity}). These results answer the corresponding questions raised in \cite[Section 5]{maclagan2022higher}.

The paper is organized as follows. \Cref{sec:balanced} introduces balanced top cycles and their behavior under subdivisions and quotients. \Cref{sec:filling} establishes the bound on surviving connected components using filling arguments. \Cref{sec:skeleta} proves the connectivity of lower-dimensional skeleta and \Cref{thm:general}. \Cref{sec:trop-lin} applies the results to tropical linear spaces and Bergman fans. \Cref{sec:algebraic} examines coordinate-projection matroids and related examples.

\section{Balanced cycles and polyhedral structures}
\label{sec:balanced}

Throughout, $V$ is a finite-dimensional real vector space. A polyhedral complex is a finite collection of closed convex polyhedra satisfying the usual face and intersection conditions. We make no rationality assumption; this allows nonrational subdivisions and quotients by arbitrary real subspaces. For a polyhedral complex $\Delta$, its \emph{support} is $|\Delta|=\bigcup_{\sigma\in\Delta}\sigma$, and $\Delta_q$ denotes the set of its $q$-dimensional cells. If $\Delta$ is \emph{pure} of dimension $d$, its \emph{facets} and \emph{ridges} are its cells of dimension $d$ and $d-1$, respectively. Its \emph{$k$-skeleton} consists of the cells of dimension at most $k$; if $\Delta$ is pointed, it is pure of dimension $k$. The \emph{lineality space} of a polyhedron $\sigma$ is the largest linear subspace $U$ with $\sigma+U=\sigma$. A polyhedron is \emph{pointed} if its lineality space is zero, equivalently if it contains no affine line, and a complex is pointed if all its cells are. If a linear subspace $U$ lies in the lineality space of every cell of $\Delta$, the images of the cells in $V/U$ form a polyhedral complex $\Delta/U$ with the same face poset.

\begin{defi}\label{def:connectivity}
A hypergraph is \emph{connected} if any two vertices are joined by a sequence of vertices in which consecutive vertices lie in a common hyperedge; we regard the empty hypergraph as connected. A hypergraph is \emph{strongly $p$-vertex-connected} if it remains connected after deleting any at most $p-1$ vertices together with every hyperedge containing one of them. The \emph{facet-ridge hypergraph} of a pure polyhedral complex $\Delta$ has the facets of $\Delta$ as vertices and, for each ridge $\tau$, the hyperedge consisting of the facets containing $\tau$. We say that $\Delta$ is \emph{$k$-connected through codimension one} if its facet-ridge hypergraph is strongly $k$-vertex-connected.
\end{defi}

Thus, deleting a facet removes every ridge it contains: a ridge of a deleted facet cannot be used to connect surviving facets, even if two other facets containing it survive. We call this \emph{strong deletion}. 

Let $\Delta$ be a pure $d$-dimensional polyhedral complex in $V$. For each cell $\sigma\in\Delta$, let $T_\sigma=\spanz(\sigma-\sigma)$ denote its \emph{direction space}. Fix an orientation for each cell of $\Delta$. When $V=N\otimes\RR$ for a lattice $N$, we call $T_\sigma$ \emph{rational} if it is spanned by vectors in $N$. For a fixed real vector space $A$, the \emph{cellular Borel-Moore chains} of $\Delta$ with coefficients in $A$ are
\[C_q^{\mathrm{BM}}(\Delta;A)=\bigoplus_{\sigma\in\Delta_q}A[\sigma]\]
with boundary operator
\[\partial(a[\sigma])=\sum_{\tau\prec\sigma}\epsilon_{\sigma\tau}a[\tau]\]
where the sum is over codimension-one faces $\tau$ of $\sigma$ and $\epsilon_{\sigma\tau}$ is the incidence sign. 

\begin{rmk}
Unlike ordinary cellular chains, Borel-Moore chains allow unbounded cells. We use the same coefficient space $A$ for every cell, so the boundary maps involve only incidence signs, unlike tropical chains whose coefficient spaces depend on the cells \cite{aksnes2023tropical}.
\end{rmk}

Taking $A=\bigwedge^dV$, we restrict the boundary map to the subspace $\bigoplus_{\sigma\in\Delta_d}\bigwedge^d T_\sigma[\sigma]$ and write
\begin{equation}\label{eq:boundary-map}
    \partial_d^\Delta:\bigoplus_{\sigma\in\Delta_d}\bigwedge^d T_\sigma[\sigma]\to C_{d-1}^{\mathrm{BM}}(\Delta;A).
\end{equation}
When it is clear from the context, we suppress the $\Delta$ and denote it by $\partial_d$. We call $\ker\partial_d^\Delta$ the \emph{top-cycle space} of $\Delta$. For each $\sigma\in\Delta_d$, choose a generator $\xi_\sigma\in\bigwedge^dT_\sigma$ compatible with the chosen orientation of $\sigma$. A cycle $z=\sum_{\sigma\in\Delta_d}w_\sigma\xi_\sigma[\sigma]\in\ker\partial_d$ is called \emph{positive} if $w_\sigma>0$ for every $\sigma\in\Delta_d$. The positivity of a cycle is independent of the chosen orientations and positive rescalings of $\xi_\sigma$. We call $\Delta$ \emph{uniquely balanced} if it admits a positive top cycle and $\dim\ker\partial_d^\Delta=1$. For rational fans, this agrees with \cite[Definition 3.12]{aksnes2023tropical} since in top degree the tropical Borel-Moore homology there coincides with $\ker\partial_d$.

Suppose that $V=N\otimes\RR$ for a lattice $N$ and that all direction spaces are rational. Choose primitive lattice generators $\xi_\sigma$ of $\bigwedge^d(T_\sigma\cap N)$ compatible with the orientations. Then, $\partial z=0$ is the usual balancing condition $\sum_{\tau\subset\sigma}w_\sigma u_{\sigma/\tau}=0$ in $V/T_\tau$ at each ridge $\tau$ where $u_{\sigma/\tau}$ is the primitive normal vector of $\tau$ in $\sigma$ \cite[Definition 3.12 and Example 3.14]{aksnes2023tropical}.

\begin{defi}
A positive balanced polyhedral support with rational direction spaces is \emph{tropically irreducible} if it cannot be expressed as the union of two proper positive balanced supports of the same dimension.
\end{defi}

To relate tropical irreducibility to the uniquely balanced condition, we show that the top-cycle space is preserved under subdivision and quotient.

\begin{lem}\label{lem:cycleStructure}
Let $\Delta$ be a pure $d$-dimensional polyhedral complex and $\Delta'$ a subdivision of $\Delta$. Then, there is a natural isomorphism $\ker\partial_d^\Delta\cong\ker\partial_d^{\Delta'}$ preserving positivity. Moreover, if a real subspace $U$ lies in the lineality space of every cell, quotienting gives an isomorphism $\ker\partial_d^\Delta\cong\ker\partial_{d-\dim U}^{\Delta/U}$ preserving positivity.
\end{lem}
\begin{proof}
Suppose that a maximal cell $\sigma$ is subdivided into $\sigma_1,\dots, \sigma_m$ with inherited orientations and the same exterior generator $\xi_{\sigma_i}=\xi_\sigma$. This is possible since $T_{\sigma_i}\subset T_\sigma$ for all $i$. Let $w_\sigma\xi_\sigma$ be the coefficient of $\sigma$ in the cycle. For each $\sigma_i$, we impose the same exterior coefficients $w_{\sigma_i}\xi_{\sigma_i}=w_\sigma\xi_\sigma$ for all $i=1,\dots, m$. Let $\tau$ be a new interior ridge shared by two subfacets $\sigma_i$ and $\sigma_j$. Then, we have $\epsilon_{\sigma_i\tau}=-\epsilon_{\sigma_j\tau}$, and hence $\epsilon_{\sigma_i\tau}w_\sigma\xi_\sigma+\epsilon_{\sigma_j\tau}w_\sigma\xi_\sigma=0$. Along the subdivided original ridges, the boundary equations agree with the original ones and hence vanish.

Conversely, let $z'=w_1\xi_\sigma[\sigma_1]+w_2\xi_\sigma[\sigma_2]+\cdots$ be a cycle on $\Delta'$ with coefficient $w_i\xi_{\sigma}$ on the subfacets of $\sigma$. At any new interior ridge $\tau$ shared by $\sigma_i$ and $\sigma_j$, the cycle condition gives $0=(\partial z')|_\tau=\epsilon_{\sigma_i\tau}(w_i-w_j)\xi_\sigma$, and hence $w_i=w_j$. Since the subfacets of $\sigma$ are connected through interior ridges, we have $w_1=\cdots=w_m$. Thus, a subdivision does not induce an independent weight, and the cycle spaces are naturally isomorphic. Since the weights are unchanged, positivity is also preserved.

For the quotient, fix an oriented volume element $\eta_U\in\bigwedge^{\dim U}U$. The map
\[\bar v_1\wedge\cdots\wedge\bar v_{d-\dim U}\mapsto\eta_U\wedge v_1\wedge\cdots\wedge v_{d-\dim U}\]
identifies the top exterior lines. Note that it is independent of the lift. With product orientations, the boundary equations correspond up to a common sign, so zero boundaries and positivity are preserved. Since every cell is invariant under $U$, quotienting preserves the face and incidence relations. 
\end{proof}

\begin{lem}\label{lem:irreducible}
For a complex $\Delta$ with rational direction spaces admitting a positive balanced rational weighting, tropical irreducibility is equivalent to the space of top-dimensional cycles being one-dimensional, that is, $\dim\ker \partial_d=1$.
\end{lem}
\begin{proof}
Since the boundary map $\partial_d$ is represented by a rational matrix, its kernel has a basis over $\mathbb Q$. Let $w$ be a positive rational weighting and $a$ a rational balanced weighting not proportional to $w$. Let $s=\min_\sigma \frac{a_\sigma}{w_\sigma}$ and $t=\max_\sigma \frac{a_\sigma}{w_\sigma}$. Then, we have $s<t$ and $a-sw$ and $tw-a$ are nonnegative balanced weightings. Each has nonempty proper support and their supports cover $|\Delta|$. A contradiction to tropical irreducibility. 

Conversely, a proper balanced sub-support gives a nonproportional weighting after a common refinement; the subdivision assertion in \Cref{lem:cycleStructure} identifies it with a weighting on $\Delta$, and hence $\dim\ker\partial_d^\Delta\geq2$. This is also the equivalence in \cite[Proposition 2.4 (i) and (iii)]{babaee2024extremal}.
\end{proof}

For rational direction spaces, the uniquely balanced condition is also equivalent to strong extremality of balanced weighted complexes in the sense of \cite[Definition 2.11]{babaee2015tropical}.

\section{Filling cycles on deleted facets}\label{sec:filling}

Our goal is to bound the number of connected components remaining after deleting facets by the dimension of the top-cycle space. The idea is to construct an independent cycle from each surviving component. Restricting a positive cycle to a surviving component may destroy the balancing condition, but the resulting boundary is supported on the union of the deleted facets. We first show that such boundaries can be ``filled''. The bound of $d-1$ deleted facets arises from two complementary constraints: one topological and the other geometric.

\begin{lem}\label{lem:BMunion}
Let $P_1,\dots,P_m$ be pointed closed convex polyhedra, $m\geq1$. For any finite-dimensional real vector space $A$,
\[H_q^{\mathrm{BM}}\left(\bigcup_{i=1}^mP_i;A\right)=0\]
for $q\geq m$.
\end{lem}
\begin{proof}
We first consider a single convex polyhedron. A bounded nonempty convex polyhedron is contractible, so its Borel-Moore homology vanishes in positive degrees. We show that an unbounded pointed polyhedron $P$ has zero Borel-Moore homology in every degree. To show this, choose $h$ strictly positive on $\rec(P)\setminus\{0\}$, and $R$ larger than its values at all vertices. The hyperplane $\{h=R\}$ cuts off the unbounded parts without removing any vertices. Then, $P\cap\{h\leq R\}$ and $P\cap\{h=R\}$ are nonempty compact convex polyhedra. The cellular Borel-Moore chain complex of $P$ is isomorphic to the relative cellular chain complex of $(P\cap\{h\leq R\},P\cap\{h=R\})$ since the faces of $P$ correspond to the truncated faces not contained in $P\cap\{h=R\}$. Because both $P\cap\{h\leq R\}$ and $P\cap\{h=R\}$ are contractible, their relative homology vanishes. The same argument applies to any finite subdivision by choosing $R$ above all its vertices.

We now induct on $m$. Let $D=P_1\cup\cdots\cup P_{m-1}$. The cellular Mayer-Vietoris sequence for $D\cup P_m$ contains
\[H_q^{\mathrm{BM}}(D;A)\oplus H_q^{\mathrm{BM}}(P_m;A)\rightarrow H_q^{\mathrm{BM}}(D\cup P_m;A)\rightarrow H_{q-1}^{\mathrm{BM}}(D\cap P_m;A).\]
For $q\geq m$, both $H_q^{\mathrm{BM}}(D;A)$ and $H_q^{\mathrm{BM}}(P_m;A)$ vanish by the induction hypothesis and the base case, respectively. Moreover, $D\cap P_m=\bigcup_{i=1}^{m-1}(P_i\cap P_m)$ is a union of at most $m-1$ pointed convex polyhedra, after discarding empty intersections. Thus, its $(q-1)$-th Borel-Moore homology vanishes by induction. The result follows.
\end{proof}

In our situation, if $D$ is the union of at most $d-1$ deleted $d$-dimensional facets, then $H_{d-1}^{\mathrm{BM}}(D;A)=0$ by \Cref{lem:BMunion}. Hence, for any $(d-1)$-cycle $b$ on $D$ with coefficients in $A$, there exists a $d$-chain $y$ supported on $D$ such that $\partial y=b$. However, when $A=\bigwedge^d V$, such a filling need not have coefficients in the lines $\bigwedge^dT_\sigma$ required by \eqref{eq:boundary-map}. The following lemma rules this out.

\begin{lem}\label{lem:exposed}
Let $\Delta$ be a pure pointed $d$-dimensional polyhedral complex with $d\geq2$, and let $\mathcal G$ consist of at most $d-1$ facets. Each $\sigma\in\mathcal G$ has two ridges with distinct direction hyperplanes, neither contained in another member of $\mathcal G$.
\end{lem}
\begin{proof}
Fix $\sigma\in\mathcal G$. At a vertex of $\sigma$, which exists as $\sigma$ is pointed, choose $d$ incident ridges whose normal covectors in $T_\sigma^\star$ are linearly independent. Another deleted facet can contain at most one of these ridges since its intersection with $\sigma$ is a proper face of $\sigma$.

There are at most $d-2$ other deleted facets. Hence, at least two of the selected ridges are not contained in any other member of $\mathcal G$. Their direction hyperplanes are distinct by construction.
\end{proof}

These two lemmas yield the following bound.

\begin{prop}\label{prop:components}
Let $\Delta$ be a finite pure pointed $d$-dimensional real polyhedral complex admitting a positive top cycle. After strong deletion of at most $d-1$ facets, the number of nonempty components of its surviving facet-ridge hypergraph is at most $\dim\ker \partial_d$.
\end{prop}

\begin{proof}
Fix a positive top cycle $z=\sum_{\sigma\in\Delta_d}w_\sigma\xi_\sigma[\sigma]$ with $w_\sigma>0$. If no facets are deleted, the restriction of $z$ to each facet-ridge component is already a cycle. These cycles are linearly independent, proving the claim. This also covers $d=1$.

Let $\mathcal G$ be a nonempty collection of at most $d-1$ deleted facets, and let $D$ be their union, including all faces. For each surviving component $\mathcal C$, define the chain $z_{\mathcal C}=\sum_{\sigma\in\mathcal C}w_\sigma\xi_\sigma[\sigma]$. This may not be a cycle because restricting $z$ to $\mathcal C$ may destroy the balancing equations.

However, at every ridge outside $D$, all incident facets survive and belong to the same component. Thus, the original balancing equation still holds there. Hence, $b_{\mathcal C}:=\partial z_{\mathcal C}$ is supported on $D$. Since $\partial b_{\mathcal C}=0$, it is a $(d-1)$-cycle on $D$ with constant coefficients $A=\bigwedge^dV$.

By \Cref{lem:BMunion}, $H_{d-1}^{\mathrm{BM}}(D;A)=0$. Therefore, there exists a filling $y_{\mathcal C}=\sum_{\sigma\in\mathcal G}a_\sigma[\sigma]$ with $a_\sigma\in\bigwedge^dV$ such that $\partial y_{\mathcal C}=b_{\mathcal C}$.

We show that the filling coefficients satisfy $a_\sigma\in\bigwedge^dT_\sigma$. Fix $\sigma\in\mathcal G$ and choose the two ridges $\tau_1,\tau_2$ as in \Cref{lem:exposed}. At either ridge, $\sigma$ is the only deleted facet contributing to $\partial y_{\mathcal C}$. On the other hand, every contribution to $b_{\mathcal C}$ at that ridge belongs to $(\bigwedge^{d-1}T_{\tau_i})\wedge V$. Since $\partial y_{\mathcal C}=b_{\mathcal C}$, we obtain
\[a_\sigma\in\left((\bigwedge^{d-1}T_{\tau_1})\wedge V\right)\cap\left((\bigwedge^{d-1}T_{\tau_2})\wedge V\right).\]
The two direction hyperplanes are distinct in $T_\sigma$. Choosing a basis adapted to both hyperplanes shows that
\[\left((\bigwedge^{d-1}T_{\tau_1})\wedge V\right)\cap\left((\bigwedge^{d-1}T_{\tau_2})\wedge V\right)=\bigwedge^dT_\sigma.\]
Thus, every coefficient of $y_{\mathcal C}$ belongs to the required top exterior line.

Hence, $\tilde z_{\mathcal C}=z_{\mathcal C}-y_{\mathcal C}$ is a top cycle since $\partial\tilde z_{\mathcal C}=b_{\mathcal C}-b_{\mathcal C}=0$. Since $y_{\mathcal C}$ is supported only on deleted facets, $\tilde z_{\mathcal C}$ agrees with the original positive cycle on $\mathcal C$ and vanishes on every other surviving component. If $\sum_{\mathcal C}\lambda_{\mathcal C}\tilde z_{\mathcal C}=0$, then comparison of coefficients on any facet $\sigma\in\mathcal C$ gives $\lambda_{\mathcal C}w_\sigma=0$. Since $w_\sigma>0$, we have $\lambda_{\mathcal C}=0$ for every $\mathcal C$. Hence, these cycles are linearly independent. We conclude the proof.
\end{proof}

The following example shows that the bound in \Cref{prop:components} is sharp and illustrates how the boundary created by restricting a cycle to a surviving component can be filled using the deleted facet.

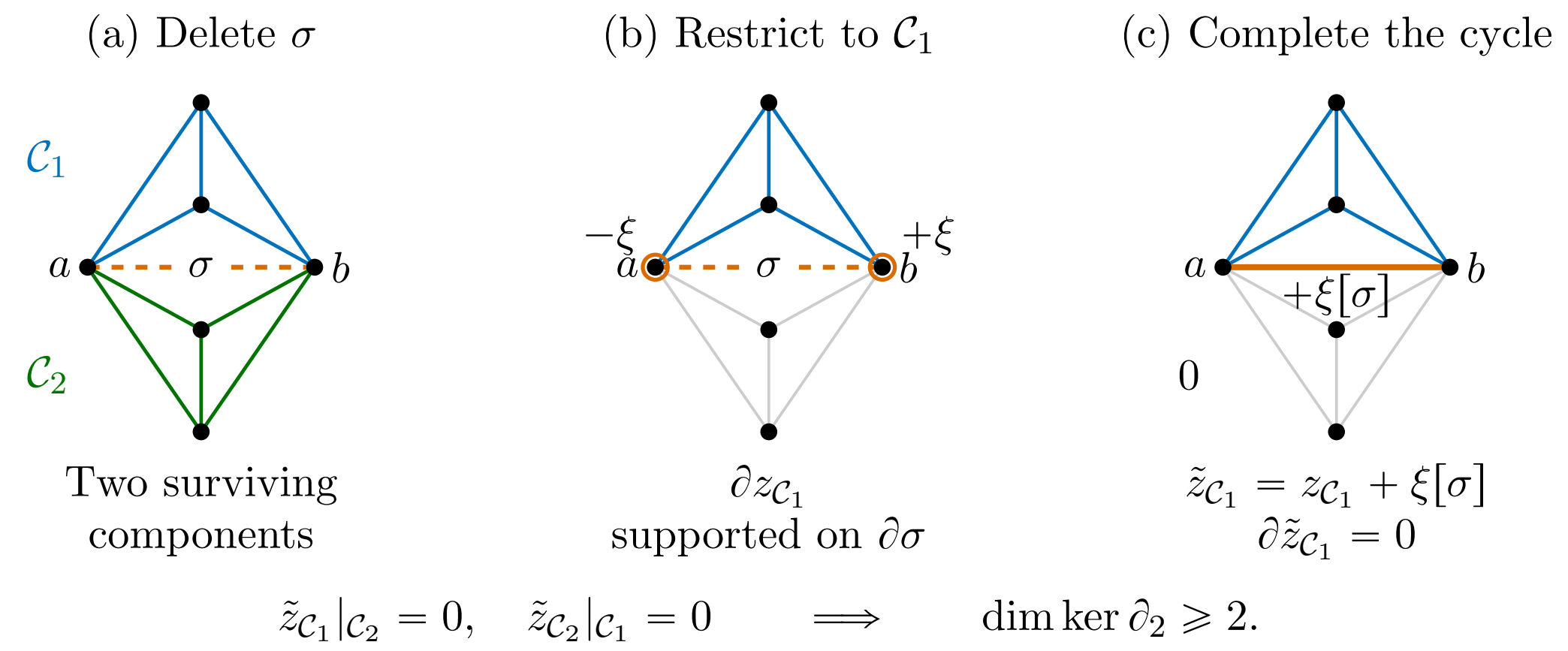
\begin{figure}[ht]
\centering
\begin{tikzpicture}[x=1cm,y=1cm,font=\small,
  componentA/.style={draw=RoyalBlue,line width=.9pt},
  componentB/.style={draw=green!45!black,line width=.9pt},
  inactive/.style={draw=black!20,line width=.7pt},
  ray/.style={circle,fill=black,inner sep=1.5pt}]
\foreach \panel/\shift in {1/0,2/5,3/10}{
  \begin{scope}[xshift=\shift cm]
    \coordinate (a) at (-1,0);
    \coordinate (b) at (1,0);
    \coordinate (c) at (0,1.45);
    \coordinate (d) at (0,.55);
    \coordinate (e) at (0,-1.45);
    \coordinate (f) at (0,-.55);
    \ifnum\panel=1
      \draw[componentA] (a)--(c)--(b)--(d)--(a) (c)--(d);
      \draw[componentB] (a)--(e)--(b)--(f)--(a) (e)--(f);
      \draw[orange!85!black,dashed,line width=1.2pt] (a)--(b);
      \node[above] at (0,1.7) {(a) Delete $\sigma$};
      \node[text=RoyalBlue] at (-1.35,.95) {$\mathcal C_1$};
      \node[text=green!45!black] at (-1.35,-.95) {$\mathcal C_2$};
      \node[fill=white,inner sep=2pt] at (0,0) {$\sigma$};
      \node[align=center] at (0,-2.15) {Two surviving\\components};
    \else
      \draw[componentA] (a)--(c)--(b)--(d)--(a) (c)--(d);
      \draw[inactive] (a)--(e)--(b)--(f)--(a) (e)--(f);
      \ifnum\panel=2
        \draw[orange!85!black,dashed,line width=1.2pt] (a)--(b);
        \draw[orange!85!black,line width=1pt] (a) circle (.11) (b) circle (.11);
        \node[above] at (0,1.7) {(b) Restrict to $\mathcal C_1$};
        \node[fill=white,inner sep=2pt] at (0,0) {$\sigma$};
        \node at (-1.4,.28) {$-\xi$};
        \node at (1.4,.28) {$+\xi$};
        \node[align=center] at (0,-2.15) {$\partial z_{\mathcal C_1}$\\supported on $\partial\sigma$};
      \else
        \draw[orange!85!black,line width=1.5pt] (a)--(b);
        \node[above] at (0,1.7) {(c) Complete the cycle};
        \node[inner sep=2pt] at (0,-.25) {$+\xi[\sigma]$};
        \node at (-1.3,-.95) {$0$};
        \node[align=center] at (0,-2.15) {$\tilde z_{\mathcal C_1}
          =z_{\mathcal C_1}+\xi[\sigma]$\\$\partial\tilde z_{\mathcal C_1}=0$};
      \fi
    \fi
    \foreach \v in {a,b,c,d,e,f}{\node[ray] at (\v) {};}
    \node[left] at (a) {$a$};
    \node[right] at (b) {$b$};
  \end{scope}
}
\node[align=center] at (5,-3.1) {
  $\tilde z_{\mathcal C_1}|_{\mathcal C_2}=0,\quad\tilde z_{\mathcal C_2}|_{\mathcal C_1}=0\qquad\Longrightarrow\qquad\dim\ker\partial_2\geq2.$};
\end{tikzpicture}
\caption{Constructing a top cycle from a surviving component in \Cref{prop:components}. Vertices represent rays and edges represent two-dimensional facets. The dashed edge represents the deleted facet $\sigma$. (a) Deleting $\sigma$ separates the two components. (b) Restricting the original cycle to one component creates boundary coefficients $-\xi$ and $+\xi$. (c) Adding $\xi[\sigma]$ cancels this boundary and completes the cycle.}
\label{fig:componentFilling}
\end{figure}

\begin{example}\label{ex:componentFilling}
For linearly independent $u,v,w$, let $P(u,v,w)$ be the fan whose rays are generated by $u,v,w,-u-v-w$ and whose facets are the cones on every pair of these rays. Set $X=P(e_1,e_2,e_3)\cup P(e_1,2e_2,e_4),\sigma=\operatorname{pos}(e_1,e_2)$ and $\xi=e_1\wedge e_2$ (see \Cref{fig:componentFilling}). The two fans meet exactly in $\sigma$ and its faces. Each fan is a linear image of the standard tropical plane. Normalize its positive cycle $z_i$ to have exterior coefficient $\xi$ on $\sigma$, and take $z=z_1+z_2$. Thus, the coefficient of $z$ on $\sigma$ is $2\xi$. Strong deletion of $\sigma$ leaves precisely two components $\mathcal C_1$ and $\mathcal C_2$. With orientations chosen so that $\partial[\sigma]=[a]-[b]$, we have $\partial z_{\mathcal C_i}=-\xi\partial[\sigma]$ and $\tilde z_{\mathcal C_i}=z_{\mathcal C_i}+\xi[\sigma]=z_i$. The top-cycle space has dimension exactly two: the balancing equations at the two nonshared rays of each plane force its five surviving facet coefficients to be proportional to those of $z_i$. After subtracting these two multiples, any remaining cycle is supported on $\sigma$ alone and hence is zero. Thus, the bound in \Cref{prop:components} is attained.
\end{example}

\Cref{lem:BMunion} also gives a lower bound on the number of facets supporting a nonzero top cycle.

\begin{lem}\label{lem:cycleCount}
A pointed pure $d$-dimensional complex admitting a nonzero top cycle has at least $d+1$ facets.
\end{lem}
\begin{proof}
If there were $m\leq d$ facets, \Cref{lem:BMunion} would give $H_d^{\mathrm{BM}}(|\Delta|;\bigwedge^dV)=0$. But a nonzero top cycle represents a nonzero homology class since there are no $(d+1)$-cells.
\end{proof}

\section{To the skeleta}\label{sec:skeleta}

We show how connectivity through codimension one descends to lower-dimensional skeleta. Related higher-connectivity results for skeleta of arbitrary convex polytopes were studied in \cite{athanasiadis2009graph,hathcock2023hypergraph,sallee1967incidence}. The main difficulty is that deleting a lower-dimensional cell may remove either a vertex or a hyperedge from the facet-ridge hypergraph of a higher-dimensional cell. The following lemma shows that either deletion drops strong vertex connectivity by at most one.

\begin{lem}\label{lem:mixed}
Let $H$ be a strongly $p$-vertex-connected hypergraph with at least $p+1$ vertices for $p\geq 2$. Then, $H-h$ is strongly $(p-1)$-vertex-connected for every hyperedge $h$. Consequently, deleting $a$ vertices strongly and then $b$ hyperedges leaves $H$ connected whenever $a+b\leq p-1$.
\end{lem}
\begin{proof}
Fix a hyperedge $h$ of $H$. First, we delete at most $p-2$ vertices and call the resulting hypergraph $K$. $K$ has at least three vertices and remains connected even after one additional strong vertex deletion. If $h$ was removed by the vertex deletions, then $K-h=K$ which is connected. Otherwise, if $K-h$ were disconnected, $h$ would meet every component because $K$ is connected. Choose $v\in h$ so that at least two components of $K-h$ remain nonempty after removing $v$: when there are two components, choose $v$ in one with at least two vertices, and when there are at least three components, any choice works. Strong deletion of $v$ also deletes $h$, contradicting connectivity of $K-v=(K-h)-v$. Iterating this observation after the first $a$ vertex deletions proves the claim.
\end{proof}

We next apply this observation to the facet-ridge graphs of pointed polyhedra.

\begin{lem}\label{lem:boundary}
Let $Q$ be a pointed polyhedron of dimension $q\geq2$. Its facet-ridge hypergraph is $(q-1)$-connected and has at least $q$ vertices.
\end{lem}
\begin{proof}
Homogenize $Q$ to a pointed cone in one higher dimension and cut by a strictly positive affine section. The resulting compact $q$-polytope $\overline Q$ has a distinguished face $F_\infty$ at height zero, possibly empty. Faces of $\overline{Q}$ not contained in $F_\infty$ correspond to the faces of $Q$ by projective transformation.

The facet-ridge graph of $\overline Q$ is the vertex-edge graph of its polar, hence is $q$-connected by Balinski's theorem \cite{balinski1961graph}. The ridges of $\overline{Q}$ contained in $F_\infty$ are not faces of $Q$, and every other ridge lies in exactly two facets of $Q$, so the facet-ridge hypergraph of $Q$ is a graph. Hence, if $\dim F_\infty=q-1$, passing to $Q$ deletes one vertex and its incident edges from the facet-ridge graph of $\overline{Q}$. If $\dim F_\infty=q-2$, passing to $Q$ deletes one edge. If its dimension is smaller, or $Q$ is bounded, it deletes neither vertices nor edges. \Cref{lem:mixed} proves $(q-1)$-connectivity. A $q$-polytope has at least $q+1$ facets, and at most one is lost, proving the second claim.
\end{proof}

We can descend from the connectivity of a polyhedral complex to that of its lower-dimensional skeleta.

\begin{figure}[ht]
\centering
\begin{tikzpicture}[
  x=1cm,y=1cm,
  font=\small,
  upper/.style={
    circle,draw=black!65,fill=white,
    minimum size=9mm,inner sep=1pt
  },
  lower/.style={
    circle,draw=RoyalBlue,fill=RoyalBlue!10,
    inner sep=1.5pt
  },
  localpath/.style={
    RoyalBlue,very thick
  },
  shared/.style={
    orange!80!black,densely dashed
  }
]

\foreach \a/\b/\i in {0/3/0,3/6/1,6/9/2}{
  \fill[RoyalBlue!5]
    (\a,-.5) rectangle (\b,.45);
  \draw[black!20,dashed]
    (\a,-.5) rectangle (\b,.45);
  \node at ({(\a+\b)/2},-.82)
    {$\text{inside }Q_{\i}$};
}

\node[above] at (4.5,3.25)
  {Global connectivity of $(k+1)$-cells};

\node[upper] (Q0) at (1.5,2.3) {$Q_0$};
\node[upper] (Q1) at (4.5,2.3) {$Q_1$};
\node[upper] (Q2) at (7.5,2.3) {$Q_2$};

\draw[->,thick,black!65]
  (Q0) -- node[above] {$\rho_1$} (Q1);
\draw[->,thick,black!65]
  (Q1) -- node[above] {$\rho_2$} (Q2);

\node[above] at (4.5,.65)
  {Connectivity of surviving $k$-cells};

\node[lower,label=below:{$\alpha$}] (a)
  at (0,0) {};
\node[lower] (u0) at (1.5,0) {};
\node[lower,label=below:{$\rho_1$}] (g1)
  at (3,0) {};
\node[lower] (u1) at (4.5,0) {};
\node[lower,label=below:{$\rho_2$}] (g2)
  at (6,0) {};
\node[lower] (u2) at (7.5,0) {};
\node[lower,label=below:{$\beta$}] (b)
  at (9,0) {};

\draw[localpath]
  (a)--(u0)--(g1)
  (g1)--(u1)--(g2)
  (g2)--(u2)--(b);

\draw[shared,thick] (3,2.3)--(g1);
\draw[shared,thick] (6,2.3)--(g2);


\end{tikzpicture}
\caption{
Schematic illustration of the proof of \Cref{prop:descent}. Connectivity of the $(k+1)$-skeleton gives a path $Q_0,Q_1,Q_2$ where consecutive cells share surviving $k$-faces $\rho_1,\rho_2$. Within each $Q_i$, local facet-ridge connectivity provides a path between the relevant $k$-faces. Concatenating these paths connects $\alpha$ to $\beta$ in the $k$-skeleton.
}
\label{fig:descent}
\end{figure}
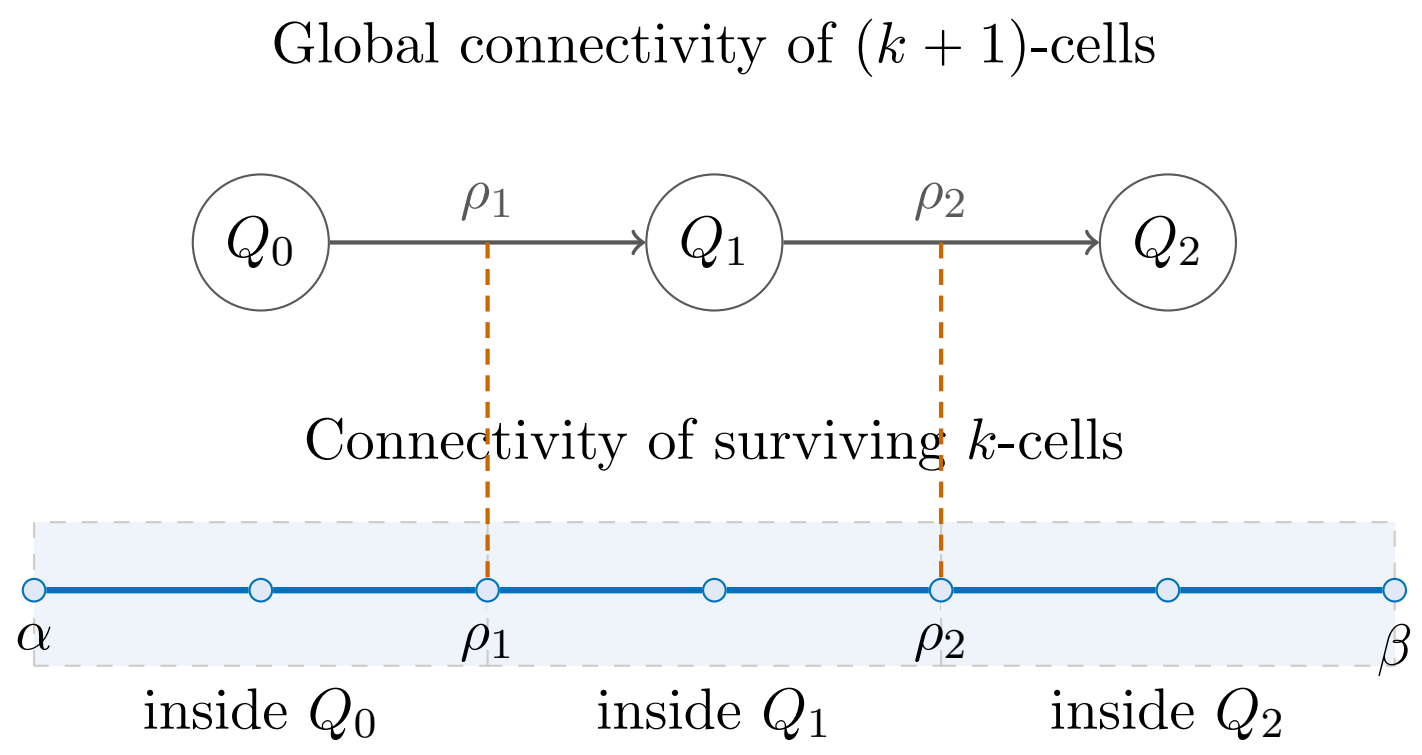

\begin{prop}\label{prop:descent}
Let $\Delta$ be a finite pure pointed $d$-dimensional polyhedral complex with at least $d+1$ facets. If $\Delta$ is $d$-connected through codimension one, then its $k$-skeleton is $k$-connected through codimension one for every $1\leq k\leq d$.
\end{prop}
\begin{proof}
Use descending induction, assuming the $(k+1)$-skeleton is $(k+1)$-connected. Delete at most $k-1$ closed $k$-cells, denoted $\mathcal G$. In the facet-ridge hypergraph of the $(k+1)$-skeleton, these $k$-cells are hyperedges. Deleting them preserves connectivity by \Cref{lem:mixed}. This hypergraph has at least $k+2$ vertices: when $k=d-1$, this is assumed; for $k\leq d-2$, we apply \Cref{lem:boundary} to any $(k+2)$-cell.

To turn this connectivity of $(k+1)$-cells into connectivity of their $k$-dimensional faces, we need to connect surviving $k$-faces within each $(k+1)$-cell. Fix such a $(k+1)$-cell $Q$. By \Cref{lem:boundary}, its facet-ridge graph is $k$-connected and has at least $k+1$ vertices. For each $\gamma\in\mathcal{G}$, the intersection $\gamma\cap{Q}$ is a single proper face or is empty. If it has dimension $k$, it deletes one graph vertex and its incident edges in the facet-ridge graph. If it has dimension $k-1$, it deletes one edge. A smaller intersection deletes neither. Thus, as we deleted at most $k-1$ $k$-cells, at most $k-1$ vertex or edge deletions occur. By \Cref{lem:mixed}, the surviving $k$-faces of $Q$ remain connected through allowed $(k-1)$-faces.

Given two surviving $k$-cells $\alpha$ and $\beta$, choose $(k+1)$-cells $Q_\alpha$ and $Q_\beta$ containing $\alpha$ and $\beta$, respectively. Because $\Delta$ is pure, such cells exist. By connectivity of the $(k+1)$-skeleton after the hyperedge deletions, $Q_\alpha$ and $Q_\beta$ can be joined by a sequence of $(k+1)$-cells, with consecutive cells sharing a surviving $k$-cell (see \Cref{fig:descent}). Within each successive $(k+1)$-cell along this sequence, connect the incoming and outgoing $k$-faces using the local connectivity established above. Concatenating these paths gives a path between the original two $k$-cells. This also works for $k=1$, when there are no deletions, and proves the induction.
\end{proof}

We now combine the connectivity bound of \Cref{prop:components} with the descent result of \Cref{prop:descent}. When the top-cycle space is one-dimensional, the resulting connectivity extends to every skeleton after quotienting by the common lineality space.

\begin{proof}[Proof of \Cref{thm:general}]
Restricting a positive top cycle to each facet-ridge component gives a nonzero top cycle since all facets incident to a ridge belong to the same component. These cycles have disjoint supports, and so they are linearly independent. Since $\dim\ker\partial_d^\Delta=1$, there can be only one component. A nonempty face of a polyhedron has the same lineality as that polyhedron. Hence, facets sharing a ridge have the same lineality space. Then, connectivity through codimension one and purity imply that all cells of $\Delta$ share a common lineality space $U$. By \Cref{lem:cycleStructure}, the pointed quotient has a positive top cycle and a one-dimensional top-cycle space. Apply \Cref{prop:components,lem:cycleCount,prop:descent} in dimension $d-\dim U$. When $d-\dim U=0$, the connectivity assertion is vacuous.
\end{proof}

For complexes with rational direction spaces, unique balancing is equivalent to tropical irreducibility by \Cref{lem:irreducible}. Thus, we obtain the following consequence.

\begin{cor}\label{cor:irreducible}
A finite pure tropically irreducible balanced complex of dimension $d$, with common lineality of dimension $\ell$, is $(d-\ell)$-connected through codimension one. Its $(\ell+k)$-skeleton is $k$-connected through codimension one for $1\leq k\leq d-\ell$.
\end{cor}
\begin{proof}
Apply \Cref{lem:irreducible,thm:general}.
\end{proof}

The general theorem also applies to arbitrary pointed polyhedral subdivisions of $\RR^n$.

\begin{cor}\label{cor:euclideanSubdivision}
Let $\Delta$ be a finite pointed polyhedral complex with support $\RR^n$. For every $1\leq k\leq n$, its $k$-skeleton is $k$-connected through codimension one.
\end{cor}
\begin{proof}
The complex $\Delta$ is pure of dimension $n$, and every facet has direction space $\RR^n$. Hence, $\ker\partial_n^\Delta=H_n^{\mathrm{BM}}(\Delta;\bigwedge^n\RR^n)\cong\RR$. The fundamental cycle induced by an orientation of $\RR^n$ is positive. Since $\Delta$ is pointed, \Cref{thm:general} applies with $U=0$.
\end{proof}

Applied to the normal fan of a $d$-polytope $P$, whose $(d-k)$-skeleton has the $k$-faces of $P$ as facets and the $(k+1)$-faces as ridges, \Cref{cor:euclideanSubdivision} recovers \cite[Theorem 1.1]{hathcock2023hypergraph}, and hence \cite[Corollary 2]{maclagan2022higher}. It also covers non-polytopal complete fans and pointed subdivisions of $\RR^n$ that are not fans.

\section{Tropical linear spaces and Bergman fans}\label{sec:trop-lin}

\Cref{thm:general} requires the global condition that $\dim \ker\partial_d^\Delta=1$, involving balancing equations across all facets and ridges. However, it can be obtained locally when the balancing relations at individual ridges are unique and the complex is connected through codimension one. The following proposition makes this observation precise.

\begin{prop}\label{prop:ridgeCriterion}
Let $\Delta$ be a finite pure $d$-dimensional complex with a positive top cycle $z=\sum_\sigma w_\sigma\xi_\sigma[\sigma]$. For each ridge $\tau$, define a linear map
\[R_\tau:\RR^{\{\sigma\in\Delta_d\mid\tau\subset\sigma\}}\to\bigwedge^d V\] with $(a_\sigma)_\sigma\mapsto\sum_{\tau\subset\sigma}\epsilon_{\sigma\tau}a_\sigma\xi_\sigma$. Thus, $\ker R_\tau$ consists of the weightings balanced at $\tau$. Suppose $\dim\ker R_\tau=1$ for every ridge. Then, $\dim\ker\partial_d^\Delta$ equals the number of facet-ridge components of $\Delta$. In particular, if $\Delta$ is connected through codimension one, then $\ker\partial_d^\Delta=\RR z$, and the conclusions of \Cref{thm:general} hold.
\end{prop}
\begin{proof}
Let $a=\sum_\sigma a_\sigma\xi_\sigma[\sigma]$ be a top cycle. At a ridge $\tau$, both $(a_\sigma)_{\tau\subset\sigma}$ and the strictly positive vector $(w_\sigma)_{\tau\subset\sigma}$ belong to $\ker R_\tau$. Hence, $\frac{a_\sigma}{w_\sigma}$ has the same value for all facets $\sigma$ containing $\tau$. These ratios are constant on each facet-ridge component. 

Conversely, restricting $z$ to any such component gives a top cycle because all facets incident to any ridge belong to the same component. These restrictions have disjoint facet supports and form a basis of $\ker \partial_d^\Delta$.
\end{proof}
 
\begin{rmk}\label{rmk:ridgeLine}
Under the assumptions of \Cref{prop:ridgeCriterion}, the existence of a positive top cycle and the condition $\dim\ker R_\tau=1$ imply that the transverse star
\[\mathrm{Star}_\Delta(\tau)/T_\tau=\bigcup_{\tau\subset \sigma}\mathbb{R}_{\geq 0}u_{\sigma/\tau}\]
is linearly isomorphic to the image of the standard tropical line in $\mathbb{R}^m/\mathbb R\boldsymbol{1}$ where $m$ denotes the number of facets containing $\tau$. More precisely, after rescaling by the positive balancing coefficients, the vectors $u_{\sigma/\tau}$ are the images of the standard basis vectors of $\mathbb R^m$ under an injective linear map from $\mathbb R^m/\mathbb R\boldsymbol{1}$ to $V/T_\tau$. Note that this linear isomorphism may not preserve the integral lattice, so the ridge condition does not imply tropical smoothness.
\end{rmk}

We apply this criterion to \emph{tropical linear spaces} of valuated matroids (see \cite[Section 4.4]{maclagan2015introduction} or \cite{doi:10.1137/080716219}, for example). Their local geometry is governed by initial matroids, which allows us to verify the ridge condition explicitly. Let $p=(p_B)$ be a valuated matroid of rank $r$ on $[n]$ with loopless underlying matroid $M$. For $w\in\mathbb R^n$, write $w(B)=\sum_{i\in B}w_i$, and let $M_w$ be the matroid whose bases maximize $w(B)-p_B$. The associated tropical linear space is
\[\Tp=\{w\in\mathbb R^n\mid M_w\text{ is loopless}\}\]
which is pure of dimension $r$. The \emph{Bergman fan} $\mathcal B(M)$ of a matroid $M$ is the tropical linear space associated with the trivial valuation. For its combinatorial structure, see \cite{ardila2006bergman,feichtner2005matroid}.

For $F\subset[n]$, write $e_F=\sum_{i\in F}e_i$. Let $c(M)$ denote the number of connected components of $M$, and set
\[L=\mathrm{span}\{e_C\mid C\text{ is a connected component of }M\}.\]
Thus, $\dim L=c(M)$.

For a polyhedral support $X\subset V$ and $x\in X$, write
\[\mathrm{Tan}_x X=\{v\in V \mid x+\varepsilon v\in X\text{ for all sufficiently small }\varepsilon>0\}.\]
The \emph{coarse subdivision} of $\Tp$ consists of the closures of the loci on which $M_w$ is constant. Its cells are defined by comparisons of the basis scores $w(B)-p_B$.

\begin{lem}\label{lem:TLSweights}
Let $p$ be a valuated matroid of rank $r$ with a loopless underlying matroid. Then, every finite pure $r$-dimensional polyhedral complex $\Delta$ with support $\Tp$ admits a positive top cycle, and $\dim\ker\partial_r^\Delta=1$.
\end{lem}
\begin{proof}
We first claim that we have $\mathrm{Tan}_w\Tp=\mathcal B(M_w)$ for every $w\in\Tp$. To see this, for any fixed $v\in\RR^n$ and sufficiently small $\varepsilon>0$, the bases maximizing
\[(w+\varepsilon v)(B)-p_B=w(B)-p_B+\varepsilon v(B)\]
are precisely those that maximize $v(B)$ among the bases of $M_w$. Thus, $M_{w+\varepsilon v}=(M_w)_v$, which gives the claim. 

Let $\sigma$ be a cell of the coarse subdivision and $w\in\relint(\sigma)$. A vector $v$ belongs to $T_\sigma=\mathrm{Span}(\sigma-\sigma)$ precisely when $v(B)=v(C)$ for all bases $B,C$ of $M_w$. These vectors are exactly those constant on each connected component of $M_w$. Hence, $T_\sigma$ is the lineality space of $\mathcal B(M_w)$ and $\dim \sigma=c(M_w)$.

Let $\tau$ be a ridge of the coarse subdivision and $w\in\relint(\tau)$. Then, $c(M_w)=r-1$, so $M_w$ has exactly one component $S$ of rank two, and all other components have rank one. Let $P_1,\dots,P_m$ be the parallel classes of $M_w|_S$. Since $S$ is connected, $m\geq3$. 

Near $w$, the facets containing $\tau$ correspond to the rays generated by $e_{P_1},\dots,e_{P_m}$ modulo $T_\tau$. These satisfy $\sum_i e_{P_i}=e_S\in T_\tau$, and any $m-1$ of them are linearly independent modulo $T_\tau$. Hence, this is their unique linear relation up to scaling. Moreover, the vectors $e_{P_i}$ give the primitive normal vectors in the quotient lattice. Consequently, choosing the primitive lattice generators $\xi_\sigma$, the cycle $\sum_\sigma\xi_\sigma[\sigma]$ with all coefficients one is balanced at every ridge, and $\dim\ker R_\tau=1$ at every ridge.

Tropical linear spaces are connected through codimension one \cite[Corollary 1.6]{anderson2024products}. By \Cref{prop:ridgeCriterion}, the top-cycle space is spanned by the weight-one cycle. Finally, \Cref{lem:cycleStructure} transfers this to any finite polyhedral complex with support $\Tp$ through a common refinement.
\end{proof}

We verify that $L$ lies in the lineality space of every cell of the coarse subdivision. Indeed, any two bases of $M$ have the same cardinality in each connected component of $M$. Thus, translation by a vector in $L$ adds the same constant to every basis score $w(B)-p_B$, leaving the initial matroid unchanged. Consequently, every cell is invariant under translation by $L$.

\begin{thm}\label{thm:main}
Let $W$ be a linear subspace of $L$ and let $\Sigma$ be any pointed polyhedral complex with support $\Tp/W$. For every $1\leq k\leq r-\dim W$, the $k$-skeleton of $\Sigma$ is $k$-connected through codimension one.
\end{thm}
\begin{proof}
By \Cref{lem:TLSweights}, the coarse subdivision of $\Tp$ admits a positive top cycle and has a one-dimensional top cycle space. Since $W\subset L$, quotienting by $W$ preserves these properties by \Cref{lem:cycleStructure}. Then, a finite common refinement of the resulting quotient complex and $\Sigma$ transfers the same properties to $\Sigma$. Since $\Sigma$ is pointed and has dimension $d=r-\dim W$, the conclusion follows from \Cref{thm:general}.
\end{proof}

\begin{cor}\label{cor:bergmanConnectivity}
Let $M$ be a loopless matroid of rank $r$ on $[n]$, and let $L$ be the lineality space of $\mathcal{B}(M)$. Let $W\subset L$ and let $\Sigma$ be a pointed fan with support $\mathcal B(M)/W$. Put $d=r-\dim W\geq1$. Then, its $k$-skeleton is $k$-connected through codimension one for every $1\leq k\leq d$. In particular, the coarse subdivision in $\RR^n/L$ and every finer fan structure have $(r-c(M))$-connected facet-ridge hypergraphs whenever $r>c(M)$. For connected $M$ of rank at least two, the coarse fan in $\RR^n/\RR\boldsymbol{1}$ is $(r-1)$-connected.
\end{cor}
\begin{proof}
Apply \Cref{thm:main} to the trivial valuation.
\end{proof}

\section{Connectivity and coordinate-projection matroids}\label{sec:algebraic}

For a balanced complex with support $X\subset\RR^n$, let $\mathcal I(X)=\{S\subset[n]\mid\dim\pi_S(X)=|S|\}$ where $\pi_S:\RR^n\to\RR^S$ is the coordinate projection. If $X$ is the tropicalization of an irreducible variety, then $\mathcal{I}(X)$ is the collection of independent sets of its algebraic matroid \cite{yu2017algebraic}. The following examples show that neither higher connectivity nor tropical irreducibility forces $\mathcal I(X)$ to be a matroid, and that $\mathcal I(X)$ need not be algebraic even when $X$ is tropically irreducible and $\mathcal{I}(X)$ is a matroid.

\begin{example}\label{ex:nonmatroid}
Let $L_{ij}=\spanz(e_i,e_j)\subset\RR^4$ and $X=L_{12}\cup L_{23}\cup L_{34}$ with each plane subdivided into its four coordinate quadrants. Since each plane is balanced, this pointed fan is balanced with all weights one. It is $2$-connected through codimension one. In fact, after one quadrant is deleted, the other three quadrants of its plane remain connected, every other plane loses at most one boundary ray and remains connected, and consecutive planes still share a ray of their common axis. However, the two-element sets in $\mathcal I(X)$ are exactly $\{1,2\}, \{2,3\}$ and $\{3,4\}$. Thus, $\{1\}$ cannot be augmented from $\{3,4\}$, and $\mathcal I(X)$ is not a matroid. This fan is not tropically irreducible because each plane carries its own balanced weighting.
\end{example}

\begin{example}\label{ex:vamos}
Let $V_8$ be the V\'amos matroid, a connected rank-four matroid on eight elements that is not algebraic over any field \cite{ingleton1975nonalgebraic}, and let $X=\mathcal B(V_8)\subset\RR^8$. By \Cref{lem:TLSweights,lem:irreducible}, $X$ is tropically irreducible, and by \Cref{cor:bergmanConnectivity}, every pointed fan structure on $X/\RR\boldsymbol{1}$ is $3$-connected through codimension one. However, $\mathcal I(X)$ is the collection of independent sets of $V_8$ \cite[Lemma 3]{yu2017algebraic}, which is not algebraic. In particular, $X$ is not the tropicalization of any irreducible variety.
\end{example}

\begin{example}\label{ex:irreducibleNonmatroid}
Let $f_1,\dots,f_8\in\ZZ^8$ be the columns of
\[M=(m_{ij})=\begin{bmatrix}
0&5&5&3&-4&0&0&0\\
5&0&-5&3&-4&0&0&0\\
5&5&0&-3&4&0&0&0\\
5&-5&5&0&0&0&0&0\\
0&0&0&0&0&5&5&5\\
0&0&0&4&3&0&-5&5\\
0&0&0&4&3&5&0&-5\\
0&0&0&4&3&-5&5&0
\end{bmatrix}.\]
Then, $MM^{\top}=75I_8$, and so we get $\sum_jm_{ij}f_j=75e_i$ and $\sum_im_{ij}e_i=f_j$. Let $X\subset\RR^8$ be the two-dimensional fan whose maximal cones are $\mathrm{cone}(e_i,f_j)$ for $m_{ij}>0$, and $\mathrm{cone}(e_i,-f_j)$, $\mathrm{cone}(-e_i,f_j)$ and $\mathrm{cone}(-e_i,-f_j)$ for $m_{ij}<0$. This adapts the construction of \cite[Section 5.4]{babaee2015tropical}, and $M$ is obtained from two copies of the matrix in \cite[Section 5.6]{babaee2015tropical} by mixing the fourth and fifth columns. Mixing these two columns forces the scale factors of the two blocks to agree while keeping coordinates $5$ and $6$ identically zero on every cone where coordinate $4$ varies. One checks that $e_i,e_k,f_j,f_l$ are linearly independent whenever $i\neq k$ and $j\neq l$. Hence, two distinct planes $\spanz(e_i,f_j)$ and $\spanz(e_k,f_l)$ meet in $\RR e_i$, in $\RR f_j$, or at the origin, and $X$ is a pointed fan with $46$ maximal cones.

Orient each maximal cone by its two listed generators and give it weight $|m_{ij}|$. At the ray $e_i$, the weighted sum of the other generators is $\sum_jm_{ij}f_j=75e_i$, and at the ray $f_j$ it is $\sum_im_{ij}e_i=f_j$. At $-e_i$ and $-f_j$, the two cones coming from each negative entry contribute opposite vectors with equal weights. Hence, these weights define a positive top cycle. Conversely, let $w$ be any top cycle. Since $f_1,\dots,f_8$ form a basis, balancing at $e_i$ forces the weights of the cones $\mathrm{cone}(e_i,\pm f_j)$ to be $|m_{ij}|$ times a scalar depending only on $i$. Similarly, balancing at $f_j$ forces the weights of the cones $\mathrm{cone}(\pm e_i,f_j)$ to be $|m_{ij}|$ times a scalar depending only on $j$. At $-e_i$, the same argument shows that the differences between the weights of $\mathrm{cone}(-e_i,f_j)$ and $\mathrm{cone}(-e_i,-f_j)$ form a multiple of the $i$-th row of $M$. This multiple is zero because the differences vanish at a positive entry, which every row has. The same holds at $-f_j$ since every column has a positive entry. Therefore, the scalar of $i$ and the scalar of $j$ agree whenever $m_{ij}\neq0$. As the bipartite graph of nonzero entries of $M$ is connected, $w$ is a multiple of the weights above. Thus, $\dim\ker\partial_2=1$, so $X$ is tropically irreducible by \Cref{lem:irreducible} and $2$-connected through codimension one by \Cref{cor:irreducible}.

Note that $\{4\}\in\mathcal{I}(X)$, and $\{5,6\}\in\mathcal{I}(X)$ since the projection of $\mathrm{cone}(e_5,f_7)$ to the coordinates $5,6$ is two-dimensional. On the other hand, if the fourth coordinate does not vanish on $\spanz(e_i,f_j)$ with $m_{ij}\neq0$, then $i=4$ or $j\leq3$. In either case, $i\leq4$ and $j\leq3$, so this plane lies in $\spanz(e_1,\dots,e_4)$. Hence, neither $\{4,5\}$ nor $\{4,6\}$ is in $\mathcal I(X)$, and $\{4\}$ cannot be augmented from $\{5,6\}$. Therefore, $\mathcal I(X)$ is not a matroid.
\end{example}

\section*{Acknowledgments}
We thank Josephine Yu for suggesting this problem when \cite{maclagan2022higher} first appeared, and Cvetelina Hill for discussions at an early stage of this project. We also thank Diane Maclagan for her comments on the first version, which led to \Cref{ex:irreducibleNonmatroid}.

\bibliographystyle{alphaabbr}
\bibliography{ref}

\newcommand{\etalchar}[1]{$^{#1}$}
\begin{thebibliography}{GHM{\etalchar{+}}21}

\bibitem[AK06]{ardila2006bergman}
F.~Ardila and C.~J. Klivans.
\newblock The {Bergman} complex of a matroid and phylogenetic trees.
\newblock {\em Journal of Combinatorial Theory, Series B}, 96(1):38--49, 2006.

\bibitem[Aks23]{aksnes2023tropical}
E.~Aksnes.
\newblock Tropical {P}oincar{\'e} duality spaces.
\newblock {\em Advances in Geometry}, 23(3):345--370, 2023.

\bibitem[And24]{anderson2024products}
N.~Anderson.
\newblock Matroid products in tropical geometry.
\newblock {\em Research in the Mathematical Sciences}, 11(2):38, 2024.

\bibitem[Ath09]{athanasiadis2009graph}
C.~A. Athanasiadis.
\newblock On the graph connectivity of skeleta of convex polytopes.
\newblock {\em Discrete \& Computational Geometry}, 42(2):155--165, 2009.

\bibitem[Bal61]{balinski1961graph}
M.~L. Balinski.
\newblock On the graph structure of convex polyhedra in {$n$}-space.
\newblock {\em Pacific Journal of Mathematics}, 11(2):431--434, 1961.

\bibitem[BDM26]{babaee2024extremal}
F.~Babaee, S.~Dewar, and J.~Maxwell.
\newblock Extremal decompositions of tropical varieties and relations with
  rigidity theory.
\newblock {\em Journal of Symbolic Computation}, 132:102461, 2026.

\bibitem[BH17]{babaee2015tropical}
F.~Babaee and J.~Huh.
\newblock {A tropical approach to a generalized Hodge conjecture for positive
  currents}.
\newblock {\em Duke Mathematical Journal}, 166(14):2749 -- 2813, 2017.

\bibitem[BJS{\etalchar{+}}07]{bogart2007computing}
T.~Bogart, A.~N. Jensen, D.~Speyer, B.~Sturmfels, and R.~R. Thomas.
\newblock Computing tropical varieties.
\newblock {\em Journal of Symbolic Computation}, 42(1-2):54--73, 2007.

\bibitem[CP13]{cartwright2013connectivity}
D.~Cartwright and S.~Payne.
\newblock Connectivity of tropicalizations.
\newblock {\em Mathematical Research Letters}, 19(5):1089--1095, 2013.

\bibitem[FS05]{feichtner2005matroid}
E.~M. Feichtner and B.~Sturmfels.
\newblock Matroid polytopes, nested sets and {Bergman} fans.
\newblock {\em Portugaliae Mathematica}, 62(4):437--468, 2005.

\bibitem[GHM{\etalchar{+}}21]{gandini2021toric}
F.~Gandini, M.~Hering, D.~Maclagan, F.~Mohammadi, J.~Rajchgot, A.~K. Wheeler,
  and J.~Yu.
\newblock Toric and tropical bertini theorems in positive characteristic.
\newblock {\em arXiv preprint arXiv:2111.13214}, 2021.

\bibitem[HY23]{hathcock2023hypergraph}
D.~Hathcock and J.~Yu.
\newblock On the hypergraph connectivity of skeleta of polytopes.
\newblock {\em Discrete \& Computational Geometry}, 69(2):593--596, 2023.

\bibitem[IM75]{ingleton1975nonalgebraic}
A.~W. Ingleton and R.~A. Main.
\newblock Non-algebraic matroids exist.
\newblock {\em Bulletin of the London Mathematical Society}, 7(2):144--146,
  1975.

\bibitem[MS15]{maclagan2015introduction}
D.~Maclagan and B.~Sturmfels.
\newblock {\em Introduction to tropical geometry}, volume 161.
\newblock American Mathematical Society, Providence, RI, 2015.

\bibitem[MY22]{maclagan2022higher}
D.~Maclagan and J.~Yu.
\newblock Higher connectivity of tropicalizations.
\newblock {\em Mathematische Annalen}, 384:775--788, 2022.

\bibitem[Sal67]{sallee1967incidence}
G.~T. Sallee.
\newblock Incidence graphs of convex polytopes.
\newblock {\em Journal of Combinatorial Theory}, 2(4):466--506, 1967.

\bibitem[Spe08]{doi:10.1137/080716219}
D.~E. Speyer.
\newblock Tropical linear spaces.
\newblock {\em SIAM Journal on Discrete Mathematics}, 22(4):1527--1558, 2008.

\bibitem[Yu17]{yu2017algebraic}
J.~Yu.
\newblock Algebraic matroids and set-theoretic realizability of tropical
  varieties.
\newblock {\em Journal of Combinatorial Theory, Series A}, 147:41--45, 2017.

\end{thebibliography}
\end{document}